\documentclass[12pt]{amsart}

\usepackage[centertags]{amsmath}

\usepackage{amsfonts}
\usepackage{amssymb}
\def\RR{\mathbb R}
\def\QQ{\mathbb Q}

\def\into{\longrightarrow}

\def\res{\hspace{-4pt} \upharpoonright \hspace{-3pt}}
\def\harp{\res}

\DeclareMathOperator{\supp}{supp}

\DeclareMathOperator{\dom}{dom}

\DeclareMathOperator{\dep}{dep}

   \DeclareMathOperator{\Th}{Th}

\title{Deciding the Chromatic Numbers of  algebraic hypergraphs}

\author{James H. Schmerl}

\date{\today}

\begin{document}

\maketitle

\begin{abstract} For each infinite cardinal $\kappa$, the set of algebraic hypergraphs having 
chromatic number no larger than $\kappa$ is decidable. \end{abstract}

 A polynomial $p(x_0,x_1, \ldots, x_{k-1})$ over the reals $\RR$ is $(k,n)$-{\bf ary} if each $x_i$ 
is an $n$-tuple of variables.  Following \cite{avoid}, we say  that a $(k,n)$-ary polynomial $p(x_0,x_1, \ldots, x_{k-1})$ is  {\bf avoidable} if the points of $\RR^n$ can be colored with countably many colors 
 such that whenever $a_0,a_1, \ldots, a_{k-1}$ $ \in \RR^n$ 
are distinct and $p(a_0,a_1, \ldots, a_{k-1}) = 0$,  
then there are  $i < j < k$ such that the points $a_i,a_j$ are   differently  colored. The polynomial is {\bf unavoidable} if it is not avoidable. The motivating examples are the $(3,n)$-ary polynomials 
$$
\|x-y\|^2 - \|y-z\|^2, {\mbox{ for }} 2 \leq n < \omega,
$$ 
which, after a series of partial results  \cite{t1}, \cite{dav72}, \cite{t3}, \cite{ek}, \cite{t2}, \cite{t4},   were shown 
\cite{t5} to be avoidable. If one is willing to ignore a small set of solutions (i.e., those consisting of collinear points), then this result can suggestively be paraphrased as ``the sets of vertices of isosceles triangles in $\RR^n$ can be avoided". 
All the avoidable polynomials were subsequently characterized in \cite{avoid}. This characterization led to  \cite[Theorem~0.2]{avoid} stating 
 that the set of  unavoidable polynomials  over the rationals $\QQ$ is computably enumerable (c.e.). It was then explicitly conjectured in  \cite[\S3]{avoid} that this set is computable. This conjecture will be   proved here (Corollary~4.8). 
 
 Avoidability involves countable colorings. The notion of avoidability was extended to uncountable colorings in \cite{alg}: for an infinite cardinal $\kappa$, 
 the $(k,n)$-ary polynomial $p(x_0,x_1, \ldots, x_{k-1})$ is  $\kappa$-{\bf avoidable} if the points of 
 $\RR^n$ can be colored using $\kappa$ colors 
 such that whenever $a_0,a_1, \ldots, a_{k-1} \in \RR^n$ 
are distinct and $p(a_0,a_1, \ldots, a_{k-1}) = 0$,  
then there are  $i < j < k$ such that the points $a_i,a_j$ are   differently  colored. 
A polynomial is $\kappa$-{\bf unavoidable} if it is not $\kappa$-avoidable. Avoidability is the same as $\aleph_0$-avoidability. The results of 
\cite{avoid} about avoidable polynomials were  extended in \cite{alg} to $\kappa$-avoidable polynomials for each infinite cardinal $\kappa$. In particular,   the $\kappa$-avoidable 
polynomials were characterized in \cite{alg},  leading to the conclusion (although unstated in \cite{alg}) that for every infinite $\kappa$, 
the set of  $\kappa$-unavoidable polynomials  over  $\QQ$ is c.e. 

 It should be pointed  out that whether or not a given polynomial is $\kappa$-avoidable may well depend on what $2^{\aleph_0}$ is. Perhaps the quintessential  example arises from a 
 result of Fox \cite[Coro.\@ 1]{fox} stating that, for each $k < \omega$ and each ordinal $\alpha$,  the $(k+3,1)$-ary polynomial
 $$
 x_0 +x_1 + x_2 + \cdots + x_k - x_{k+1} - kx_{k+2} 
 $$
  is $\aleph_\alpha$-avoidable iff $\aleph_{\alpha+k} \geq 2^{\aleph_0}$.  Another classic  
  example ({\it cf.} \cite[Prop.\@ 1.4]{alg}) is that if $ k < \omega$, then ``the sets of vertices of orthogonal $(k+1)$-simplices in $\RR^{k+1}$ is avoidable''  iff $2^{\aleph_0} \leq \aleph_k$.

  Not only will  the conjecture from \cite{avoid}  be proved here, but so will its extension  to 
$\kappa$-avoidable polynomials. It will be proved here (Corollary~4.8) that 
 if $\kappa$ is an infinite  cardinal, then the set of $\kappa$-avoidable polynomials over $\QQ$ is computable. This result will further be extended to 
 all polynomials over $\RR$:  if $\kappa$ is an infinite  cardinal, then the  set of $\kappa$-avoidable polynomials over $\RR$ is decidable.{\footnote{Being uncountable, this set cannot be computable. 
 Nevertheless, there is a natural way to extend  {\em computable} to 
  this context, resulting in what we are calling   {\em decidable}. See Definition~4.10.}}

The title and abstract of this paper do not mention avoidable polynomials, but  instead refer to the  chromatic numbers of algebraic hypergraphs. Next, we describe the 
connections between these two concepts. 

If $2 \leq k < \omega$, then $H$ is a $k$-hypergraph if $H = (V,E)$, where $V$ (the set of its {\bf vertices}) is any set and $E$ (the  set of its {\bf edges}) is a set of $k$-element subsets of $V$.{\footnote{In general $H = (V,E)$ is a {\bf hypergraph} if $E$ is a set of subsets of $V$. What are here called $k$-{\it hypergraphs} are more usually called $k$-{\it uniform} hypergraphs. All hypergraphs occurring in this paper are $k$-hypergraphs for some $k$.}}  A function 
$\varphi : V \into C$ is a {\bf coloring} of $H$. For a cardinal $\kappa$, the coloring $\varphi$ is a $\kappa$-{\bf coloring} if $|C| \leq \kappa$ and  is a {\bf proper} coloring if it is not constant on any edge. If there is a proper $\kappa$-coloring of $H$, then $H$ is $\kappa$-{\bf colorable}.  The least $\kappa$ for which $H$ is $\kappa$-colorable is its {\bf chromatic number}  $\chi(H)$.

Suppose that   $f : (\RR^n)^k \into \RR$ (for example, $f$ might be a $(k,n)$-ary polynomial). Then the {\bf zero} $k$-hypergraph of $f$ is $(\RR^n,E)$, where $E$ is the set of $k$-element subsets  $\{a_0,a_1, \ldots, a_{k-1}\}$ of $\RR^n$ such that $f(a_0,a_1, \ldots, a_{k-1}) = 0$.  We say that a  $k$-hypergraph  is {\bf algebraic} if 
it is the zero $k$-hypergraph of a $(k,n)$-ary polynomial. 
Finally, observe that  if $p(x_0,x_1, \ldots, x_{k-1})$ is a $(k,n)$-ary polynomial,
then its zero $k$-hypergraph is $\kappa$-colorable iff $p(x_0,x_1, \ldots, x_{k-1})$ is $\kappa$-avoidable.

The rest of this paper consists of 4 sections. The first, which comprises four subsections,
contains some of the preliminary definitions and results. The relevant material from \cite{alg} is 
summarized in~\S2. The main theorem and its proofs are in \S3, and some of its   consequences about algebraic hypergraphs are presented in \S4.

\bigskip


{\bf \S1.\@ Preliminaries.} This section contains some preliminary material. 

\smallskip


{\bf 1.1.\@ The Usual Stuff.} Every ordinal is the set all smaller ordinals, and every cardinal is an inital ordinal. In particular, $\omega$ is the set of finite ordinals and $\aleph_0 = \omega$. If $\alpha$ is an ordinal and $\lambda$ is a cardinal,   then we define 
$2_\alpha^\lambda$ by recursion on $\alpha$ as follows: $2_0^\lambda = \lambda$; $2_{\alpha+1}^\lambda =
2^{2^\lambda_\alpha}$; if $\alpha$ is a limit ordinal, then $2^\lambda_\alpha = \bigcup\{2^\lambda_\beta : \beta < \alpha\}$. We let $\beth_\alpha = 2^{\aleph_0}_\alpha$.  Given a cardinal $\lambda$ and $n < \omega$, we define $\lambda^{+n}$ recursively by: $\lambda^{+0} = \lambda$ and $\lambda^{+(n+1)} = ({\lambda^+})^{+n}$. Conventionally, let $\lambda^{+\infty} \geq  \kappa$ for any cardinals $\kappa$ and~$\lambda$.

Let $X$ be a set. If  $n < \omega$, then $X^n$ is the set of $n$-tuples of elements of $X$. 
Be cautioned that if $\gamma$ is an ordinal, then the notation $\gamma^n$ will never be used to denote an ordinal, but will always denote the set of $n$-tuples of smaller ordinals. 
If $x \in X^n$, then it will often be implicit that $x = \langle x_0,x_1, \ldots, x_{n-1} \rangle$. 
As usual,  ${\mathcal P}(X)$ is the set of all subsets of~$X$, and $[X]^n = \{A \in {\mathcal P}(X) : |A| = n\}$. 

Suppose that  $X$ is linearly ordered by $<$; for example, let $X$ be an ordinal or  a set of reals. If we write 
$\{x_0,x_1, \ldots, x_{k-1}\}_< \in [X]^k$, then it is to be understood that $x_0 < x_1 < \cdots < x_{k-1}$. 
We let $X^{(n)} = \{x \in X^n : x_0 < x_1 < \cdots < x_{n-1}\}$.  

\smallskip


{\bf 1.2.\@ The Erd\H{o}s-Rado Theorem.}  This subsection reviews the Erd\H{o}s-Rado Theorem and some of its variants that will be used later. 

\bigskip

{\sc Theorem 1.1:} (The Erd\H{o}s-Rado Theorem) {\em If $\lambda$ is an infinite cardinal and $r < \omega$, then $(2^\lambda_r)^+ \into (\lambda^+)^{r+1}_\lambda$ $($that is, 
if $F : [(2^\lambda_r)^+]^{r+1} \into \lambda$, then there is $X \subseteq (2^\lambda_r)^+$ such that $|X| = \lambda^+$ and 
$F$ is constant on $[X]^{r+1})$}. 

\bigskip

{\sc Corollary 1.2:} (The Polarized Erd\H{o}s-Rado Theorem) {\em If $\lambda$ is an infinite cardinal, $r < \omega$ and $F :  ((2^\lambda_r)^+)^{r+1} \into \lambda$, then there are $X_0,X_1, \ldots,X_r \subseteq (2^\lambda_r)^+$ such that $|X_0| = |X_1| = \cdots = |X_r| = \lambda$ and $F$ is constant on $X_0 \times X_1 \times \cdots \times X_r$.}

\bigskip

{\it Proof}. Define $G : [ (2^\lambda_r)^+]^{r+1} \into \lambda$ so that if 
$A = \{a_0,a_1, \ldots, a_r\}_< \in [ (2^\lambda_r)^+]^{r+1}$, then $G(A) = F(\langle a_0,a_1, \ldots, a_r \rangle)$. Theorem~1.1 assures that there is $Y \subseteq (2^\lambda_r)^+$ such that 
$|Y| = \lambda^+$ and $G$ is constant on $[Y]^{r+1}$. Let $X_0, X_1, \ldots, X_{r} \subseteq Y$ 
be such that $|X_0| = |X_1| = \cdots = |X_r| = \lambda$ and $\sup(X_i) < \min(X_{i+1})$ for all $i < r$.
Then $F$ is constant on $X_0 \times X_1 \times \cdots \times X_r$. \qed

\bigskip

Baumgartner \cite[Theorem~1]{bau} improved the Erd\H{o}s-Rado to a canonical version. If $X$ is a set that is linearly ordered by $<$, $k< \omega$, $F$ is a function on $[X]^k$ and $C \subseteq [X]^k$, then we say that $F$ is {\bf canonical} on $C$ if there is $I \subseteq k$ 
such that whenever $\{x_0,x_1, \ldots, x_{k-1}\}_<, \{y_0,y_1, \ldots, y_{k-1}\}_< \in C$, 
then $F(x) = F(y)$ iff $x_i = y_i$ for all $i \in I$. 

\bigskip

{\sc Theorem 1.3:} (The  Canonical Erd\H{o}s-Rado Theorem) {\em If $\lambda$ is an infinite cardinal, $r < \omega$ and $F$ is a function on $ [(2^\lambda_r)^+]^{r+1}$, then there is $X \subseteq (2^\lambda_r)^+$ such that $|X| = \lambda^+$ and 
$F$ is canonical on $[X]^{r+1}$}. 

\bigskip

If $I \subseteq k < \omega$ and $X = X_0 \times X_1 \times \cdots \times X_{k-1}$, then we define the equivalence relation $\sim_I$ on 
$X$ so that if $x,y \in X$, then $x \sim_I y$ iff $x_i = y_i$ for all $i \in I$. If $\approx$ 
is an equivalence relation on $X$ and $D \subseteq X$, then we say that $\approx$ 
is {\bf canonical} on $D$ if there is $I \subseteq k$ such that $\approx$ and $\sim_I $ agree on $D$. 

\bigskip

{\sc Corollary 1.4}: (The Polarized Canonical Erd\H{o}s-Rado Theorem) 
{\em If $\lambda$ is an infinite cardinal, $r < \omega$ and $\approx$ is an equivalence relation on 
$((2^\lambda_r)^+)^{r+1}$, then there are $X_0,X_1, \ldots,X_r \subseteq (2^\lambda_r)^+$ such that $|X_0| = |X_1| = \cdots = |X_r| = \lambda$ and $\approx$ is canonical on  $X_0 \times X_1 \times \cdots \times X_r$.}

\bigskip

{\it Proof}. Let $F$ be a function on $[(2^\lambda_r)^+]^{r+1}$ such that whenever
$$
\{x_0,x_1, \ldots, x_r\}_<, \{y_0,y_1, \ldots, y_r\}_< \in [(2^\lambda_r)^+]^{r+1},
$$
then 
$$
F(\{x_0, \ldots, x_r\}) = F(\{y_0, \ldots, y_r\}){\mbox{ iff }} 
\langle x_0, \ldots,x_r \rangle \approx \langle y_0, \ldots, y_r \rangle.
$$
 Apply Theorem~1.3 to get $Y\subseteq (2^\lambda_r)^+$ such that $|Y| = \lambda^+$ and $F$ is canonical on $[Y]^{r+1}$. 
Just as in the proof of Corollary~1.2, let $X_0, X_1, \ldots, X_{r} \subseteq Y$ 
be such that $|X_0| = |X_1| = \cdots = |X_r| = \lambda$ and $\sup(X_i) < \min(X_{i+1})$ for all $i < r$.
Then $\approx$ is canonical  on $X_0 \times X_1 \times \cdots \times X_r$. \qed

\bigskip

There is one more corollary  that will be useful.

\bigskip

{\sc Corollary 1.5}: {\em If $\lambda$ is an infinite cardinal and $r,m < \omega$, then there is a cardinal $\kappa <  2^\lambda_\omega$ such that if $F : \kappa^{r+1} \into \lambda^{(m)}$, then there are $X_0,X_1, \ldots, X_r  \subseteq \kappa$ such that 
$|X_0| = |X_1| = \cdots = |X_r| = \lambda$ and whenever 
$x,y \in X_0 \times X_1 \times \cdots \times X_r$ and $i < j < m$, then 
$F(x)_i \neq F(y)_j$}. 

\bigskip

{\it Proof}.  We choose $\kappa < 2^\lambda_\omega$ to be large enough for this proof to work. 
Consider any $F : \kappa^{r+1} \into \lambda^{(m)}$. 
Let $F_0,F_1, \ldots,F_{m-1} : \kappa^{r+1} \into \lambda$ be such that 
$F(x) = \langle F_0(x), F_1(x), \ldots, F_{m-1}(x) \rangle$ for all $x \in \kappa^{r+1}$.
Observe that $F_0(x) < F_1(x) < \cdots < F_{m-1}(x) < \lambda$.  Letting $\approx_i$ be the equivalence relation on $\kappa^{r+1}$ such that if $x,y \in \kappa^{r+1}$, then $x \approx_i y$ iff 
 $F_i(x) = F_i(y)$, then (Corollary~1.4) we can   assume that each $\approx_i$ is canonical on $\kappa^{r+1}$; thus, there are $N_i \subseteq 
r+1$ such that whenever $x,y \in \kappa^{r+1}$, then 
$x \sim_{N_i} y \Longleftrightarrow F_i(x) = F_i(y)$.
We can now take $\kappa = (2^\lambda_{2r+1})^+$.  

If $I= \{i_0,i_1, \ldots,i_r\}_< \subseteq 2r+2$ and $J = \{j_0,j_1, \ldots,j_r\}_< \subseteq 2r+2$, then 
we say that $I,J$ are {\bf interlaced} if whenever $\ell < k$, then $i_\ell < j_{\ell+1}$ and $j_\ell < i_{\ell+1}$.  We will use the following ad hoc notation for this proof: if $A = \{a_0,a_1, \ldots, a_{2r+1}\}_< \in [\kappa]^{2r+2}$ and $I \in [2r+2]^{r+1}$, then $A|I = \{a_i : i \in I\}$. 

Let $G$ be a function on $[\kappa]^{2r+2}$ having finite range 
such that whenever $A = \{a_0,a_1, \ldots,a_{2r+1}\}_< \in [\kappa]^{2r+2}$ and  $B =  \{b_0,b_1, \ldots, b_{2r+1}\}_< \in [\kappa]^{2r+2}$, then 
$G(A) = G(B)$ iff  whenever $I,J \in [2r+2]^{r+1}$ and $i,j < m$, then $F_i(A|I) = F_j(A|J)$ iff $F_i(B|I) = F_j(B|J)$.
 By Theorem~1.1, 
let $Y \subseteq \kappa$ be such that $|Y| = \lambda^+$ and $G$ is constant on $[Y]^{2r+2}$. 

\smallskip

{\it Claim}: {\em Suppose that $i,j < m$ are distinct, $A \in [Y]^{2k+2}$ and $I,J \subseteq [2k+2]^{k+1}$ are interlaced. 
Then $F_i(A|I) \neq F_j(A|J)$.}

\smallskip  

Since $I,J$ are interlaced and $\approx_i$ agrees with $\sim_{N_i}$  on $\kappa^{r+1}$, we get that 
$F_i(A|I) = F_i(A|J)$. Since $F_i(A|J) \neq F_j(A|J)$, we get that $F_i(A|I) \neq F_j(A|J)$,
proving the claim.

\smallskip

 Just as in the proofs of Corollaries~1.2 and~1.4, let $X_0, X_1, \ldots, X_{r} \subseteq Y$ 
be such that $|X_0| = |X_1| = \cdots = |X_r| = \lambda$ and $\sup(X_i) < \min(X_{i+1})$ for all $i < r$.
To see that these sets are as required, 
let $x,y \in X_0 \times X_1 \times \cdots \times X_r$. Then there are $A \in [Y]^{2k+2}$ and interlaced $I,J \subseteq [2k+2]^{k+1}$ such that $\{x_0,x_1, \ldots,x_k\} = A|I$ and $\{y_0,y_1, \ldots, y_k\} = A|J$. Then the claim implies that $F(x)_i = F_i(x)\neq F_j(y) = F(y)_j$ whenever $i < j < m$. \qed

\bigskip


{\bf 1.3.\@ Hypergraphs.} The definitions of a hypergraph and  some ancillary 
notions were given in the introduction. Suppose that 
$H_1 = (V_1,E_1)$ and $H_2 = (V_2,E_2)$ are hypergraphs.  Then, $H_1$ is a {\bf subhypergraph} of $H_2$ if $V_1 \subseteq V_2$ and $E_1 \subseteq E_2$. We write $H_1 \subseteq H_2$ when $H_1$ is a subhypergraph of $H_2$. If $V_1 \subseteq V_2$ and $E_1 = E_2 \cap {\mathcal P}(V_1)$, then $H_1$ is an {\bf induced} subhypergraph of $H_2$. If $V_1 = V_2$ and $E_1 \subseteq E_2$, then $H_1$ is a {\bf spanning} subhypergraph of $H_2$. If $f : V_1 \into V_2$ is an isomorphism from $H_1$ onto a subhypergraph of $H_2$, then $f$ is an {\bf embedding} of $H_1$ into $H_2$. 
If there is an embedding of $H_1$ into $H_2$, then $H_1$ is {\bf embeddable} into $H_2$ or $H_2$ {\bf embeds} $H_1$. If  
$H_1$ is embeddable into $H_2$, then $\chi(H_1) \leq \chi(H_2)$. 

\bigskip

{\bf 1.4.\@ Algebraicity/Semialgebraicity.}
 Let   \mbox{$\RR = (\RR,+,-,\times,0,1,\leq)$}  be the ordered field of the  real numbers.   By Tarski's famous theorem,  $\Th(\RR)$, the first-order theory of $\RR$,  is exactly the same as  ${\sf RCF}$, which is the theory of the class of all real closed ordered fields, thereby proving that  $\Th(\RR)$ is decidable.    
 
Let ${\mathcal L}_{OF} = \{+,-,\times,0,1,\leq\}$ be the first-order language appropriate for ordered fields. Consider a real closed ordered field \mbox{$R = (R,+,-,\times,0,1,\leq)$}.  If $D \subseteq R$, let ${\mathcal L}_{OF}(D)$ be ${\mathcal L}_{OF}$ augmented with names for the elements 
of $D$. If $n < \omega$ and $A \subseteq R^n$, then $A$ is $D$-{\bf definable} if it is definable in $ R$ by a first-order ${\mathcal L}_{OF}(D)$-formula. If $A$ is $R$-definable in $R$, then it is $R$-{\bf semialgebraic}.  If $A \subseteq R^n$ is the zero-set of a polynomial (or, equivalently, a set of polynomials) over  $R$, then $A$ is an $R$-{\bf algebraic} set.  We say that a set is semialgebraic 
(algebraic) when it is $\RR$-semialgebraic ($\RR$-algebraic).

Some definitions from the introductions are generalized from $\RR$ to an arbitrary real closed field $R$. If $f : (R^n)^k \into R$, then the {\bf zero} $k$-hypergraph of $f$ is $(R^n,E)$, where $E$ is the set of $k$-element subsets $\{a_0,a_1, \ldots,a_{k-1}\}$ of $R^n$ such that $f(a_0,a_1, \ldots, a_{k-1}) = 0$. A $k$-hypergraph is $R$-{\bf algebraic} if it is the zero hypergraph of some $(k,n)$-ary polynomial over $R$. 

If $m < \omega$, then an {\bf open} $m$-{\bf cube} is a set $B = B_0 \times B_1 \times \cdots \times B_{m-1} \subseteq \RR^m$, where each $B_i$ is a nonempty open interval of $\RR$. For example, $\RR^m$ 
is an open $m$-cube. More generally,  if $R$ is a real closed  field, then  a subset 
$(a,b) = \{x \in R : a < x < b\}$, where $a,b \in R \cup \{-\infty,\infty\}$, is  an open interval of ${ R}$. An {\bf open} $m$-{\bf cube of} 
${ R}$ is a subset   $B = B_0 \times B_1 \times \cdots \times B_{m-1} \subseteq R^m$, where each $B_i$ is a nonempty open interval of ${ R}$.

 Suppose that  $p(x_0,x_1, \ldots, x_{k-1})$ is a $(k,n)$-ary polynomial. Then 
$p(x_0,x_1, \ldots, x_{k-1})$ is {\bf symmetric} if whenever $\pi : k \into k$ is a permutation, then 
$p(a_0,a_1, \ldots, a_{k-1}) = p(a_{\pi(0)}, a_{\pi(1)},\ldots, a_{\pi(k-1)})$. 
It is {\bf  reflexive} if $p(a_0,a_1, \ldots, a_{k-1}) = 0$ whenever $a_i = a_j$ for some $i < j < k$. 
 If  $H$ is the zero hypergraph of the polynomial $p(x_0,x_1, \ldots,$ $ x_{k-1})$, then $H$ 
 is also the zero hypergraph of the polynomial 
 $$
 \prod_\pi p(x_{\pi(0)}, x_{\pi(1)},\ldots, x_{\pi(k-1)}) \times \prod_{i<j<d}(x_i - x_j),
 $$ which is both symmetric and reflexive.    
 
 If $f : U \into \RR$, where $U \subseteq \RR^m$, then $f$ is a {\bf Nash} function if $U$ is an open set and $f$ is both semialgebraic and $C^\infty$. A function  $f : U \into \RR^n$ is Nash if each of its $n$ component functions is Nash. These definitions do generalize to any real closed  field $R$, in which case we say that $f$ is $R$-{\bf Nash}.  The following lemma refers to $R$-semialgebraically 
 connected sets (see, for example, \cite[Chap.\@ 3.2]{rag}).

 \bigskip

 {\sc Lemma 1.6}: {\em Let $f : U \into \RR$ be a Nash function, where $U \subseteq \RR^n$ is connected. Suppose that $D \subseteq \RR$, $f$ is $D$-definable, and there is $t = \langle  t_0,t_1, \ldots,t_{n-1} \rangle  \in U$  such that  $f(t_0,t_1, \ldots,t_{n-1}) = 0$ 
 and $t_0,t_1, \ldots,t_{n-1}$ are algebraically independent  over $D$. Then, $f$ is identically $0$ 
 on $U$.}  
 
 \bigskip
 
 {\it Proof}. Suppose that $t \in U$ is as in the hypothesis. Consider any $a \in U$ intending to prove that $f(a) = 0$. Since $U$ is connected, there is a semialgebraic path in $U$ from $t$ to $a$ 
 (\cite[Prop.\@ 2.5.13]{rag}), and since $U$ is open there is a rectilinear path (that is, the union of 
 finitely many line segments parallel to the coordinate axes) in $U$  from $t$ to $a$ 
 such that each endpoint $s = \langle s_0,s_1, \ldots,s_{n-1} \rangle$ (except for $a$) of each segment. 
 is such that $s_0,s_1, \dots, s_{n-1}$ are algebraically independent over $D$. Then, using induction on $n$, we can assume that $n = 1$. Thus, $U \subseteq \RR$ is an open interval. Since $t_0$ is not algebraic over $D$ and $f$ is $D$-definable, there are infinitely many $b \in U$ such that $f(b) = 0$. Because $f$ is analytic,  then $f \equiv 0$ on $U$, so $f(a) = 0$.  \qed
 
 \bigskip

 We will be considering a real-closed  field ${ R} \succ { \RR}$; that is, ${ R}$ is an elementary extension of ${ \RR}$.  If $m < \omega$ and $A \subseteq \RR^m$ is  semialgebraic, 
  then we let $A^R \subseteq R^m$ be defined in ${ R}$ by the same formula that defines $A$ in ${ \RR}$. This definition of  $A^R$ does not depend on the choice of the formula defining $A$, so $A^R$ is well defined. If $H = (\RR^n,E)$ is an algebraic $k$-hypergraph, 
  then $H^R = (R^n,E^R)$ is also a $k$-hypergraph.

 \bigskip
 
 {\sc Lemma 1.6}: {\em Suppose that $R \succ { \RR}$. Let $f : U \into R$ be an $R$-Nash function, where $U \subseteq R^n$ is $R$-semialgebraically connected.  Let $D \subseteq R$ be such that $f$ is $D$-definable. Suppose that  $t_0,t_1, \ldots,t_{n-1} \in R$ are such that $\langle t_0,t_1, \ldots,t_{n-1} \rangle \in U$, $f(t_0,t_1, \ldots,t_{n-1}) = 0$ and $t_0,t_1, \ldots,t_{n-1}$ are algebraically independent  over $D$. Then, $f$ is identically $0$  
 on $U$.}  
 
 \bigskip
 
 {\it Proof.} Let $Z \subseteq U$ be the union of all open subsets of $U$ on which $f$ is constantly $0$. 
 Then $t \in U$. Let $V \subseteq Z$ be the $R$-semialgebraically connected component of $Z$ to which $t$ belongs. 
 Then $V$ is open and $R$-semialgebraic. On the other hand, $V$ is relatively closed in $U$. Hence $V = U$ since $U$ is $R$-semialgebraically connected, so  $f$ is identically~$0$.
 \qed

\bigskip
                 
                 
{\bf \S2.\@ A Summary.} This section summarizes the relevant results of \cite{alg}. It also contains the 
requisite definitions for understanding these results.

If $1 \leq d < \omega$ and $2  \leq  k < \omega$, then $P$ is a $d$-{\bf dimensional} $k$-{\bf template} if $P$ is a set of $d$-tuples and $|P| = k$. Two $d$-dimensional $k$-templates $P,Q$ are {\bf isomorphic} if there is a bijection $f : P \into Q$ such that whenever $x,y \in P$  and $i < d$, 
then $x_i = y_i$ iff  $f(x)_i = f(y)_i$. 
  If both $P$ and $Q$ are $d$-dimensional $k$-templates, then we say that $Q$ is a {\bf homomorphic image} of~$P$ if there is a surjective  function $f : P \into Q$ such that   whenever  $x,y \in P$, $i < d$ 
and $x_i = y_i$, then $f(x)_i = f(y)_i$.

If  $X = X_0 \times X_1 \times \cdots \times X_{d-1}$ and $P$ is  a $d$-dimensional $k$-template, then its {\bf template hypergraph} $L(X,P)$ on $X$ is the $k$-hypergraph whose set of vertices is $X$ and whose edges are those $k$-templates $Q \subseteq X$ that are homomorphic images of~$P$. 
If $P$ is a $d$-dimensional $k$-template, then $L(\RR^d,P)$ is an algebraic $k$-hypergraph.

Let $P$ be a $d$-dimensional $k$-template. We say that a  subset $I \subseteq d$ is a {\bf distinguisher} for 
$P$ if whenever $x,y \in P$ are distinct, then $x_i \neq y_i$ for some $i \in I$. We then define 
$e(P)$ to be the least $e$ that is the cardinality of a distinguisher. Obviously, $e(P) \leq d$ since $d$ 
itself is a distinguisher. One easily proves by induction on $k$ that $1 \leq e(P) \leq k-1$. 

Theorem~1.1. of \cite{alg} asserts that  if $P$ is a $d$-dimensional $k$-template,  then $\chi\big(L(\RR^d,P)\big)$ is the least $\kappa$ such that 
$\kappa^{+(e(P)-1)} \geq 2^{\aleph_0}$. This theorem was stated to apply only to $\RR$ since the  primary interest in  \cite[\S1]{alg}
was with $L(\RR^d,P)$. However, the following more general theorem  could 
just as easily have been inferred from results in \cite{alg}.

\bigskip

{\sc Theorem 2.1}: ({\it cf}.\@ \cite[Theorem~1.1] {alg}) {\em Suppose that  $P$ is a $d$-dimensional $k$-template and $X$ is an infinite set.  Then $\chi\big(L(X^d,P)\big)$ is the least $\kappa$ such that 
$\kappa^{+(e(P)-1)} \geq |X|$.}

\bigskip

The next easily proved lemma shows that in certain situations the only $d$-dimensional $k$-templates that need to be considered are those with $d < k$. 

\bigskip

{\sc Lemma 2.2}: ({\it cf.} \cite[Lemma~1.7]{alg}) {\em Suppose that $P$ is a $d$-dimensional $k$-template. There is an $e(P)$-dimensional $k$-template $Q$ such that for every set $X$,
$L(X^{e(P)},Q)$ is embeddable into $L(X^d,P)$. $\big($Moreover, $L(X^{e(P)},Q)$ is isomorphic to an induced subhypergraph of $L(X^d,P)\big)$.}

\bigskip

 If  $A = A_0 \times A_1 \times \cdots \times A_{d-1}$, then a function $g : A \into Y$ is 
{\bf one-to-one in each coordinate} if 
whenever $a,b \in A$ are such that $a_i \neq b_i$ for exactly one $i < d$, then 
$g(a) \neq g(b)$.
 Suppose that $P$ is a \mbox{$d$-dimensional} $k$-template, $A = A_0 \times A_1 \times \cdots \times A_{d-1}$ and $H = (V,E)$ is a $k$-hypergraph.  
A function $f : A \into V$ is an {\bf immersion} of $L(A,P)$ into $H$ 
if $f $ is  one-to-one in each coordinate and is such that 
whenever $\{x_0,x_1, \ldots, x_{k-1}\}$ is an edge of $L(A,P)$ and 
$f(x_0),f(x_1), \ldots, f(x_{k-1})$ are pairwise distinct, then $\{f(x_0),f(x_1), \ldots, f(x_{k-1})\}$ 
is an edge of $H$. If there is an immersion of $L(A,P)$ into $H$, then $L(A,P)$ is {\bf immersible} in $H$. If $R$ is a real closed  field, $B \subseteq R^k$ is an open $k$-cube 
and $H = (R^n,E)$,
then it makes sense to refer to an immersion of $L(B,P)$ into $H$ as being $R$-semialgebraic or $R$-Nash.{\footnote{Indeed, what was referred to as an {\it immersion} in \cite{alg} and \cite{avoid} is what is here is being  
called a {\it Nash immersion}}} If there is such an  $R$-semialgebraic immersion, then we say that 
$L(B,P)$ is $R$-{\bf semialgebraically immersible} into $H$. If $L(B,P)$ is  $R$-semialgebraically  immersible into $H$, then there is an $R$-Nash immersion  of $L(\RR^d,P)$ into $H$.

\bigskip

{\sc Lemma 2.3:}  (\cite[Lemma~2.1]{alg}) {\em If $H$ is an algebraic $k$-hypergraph, $P$ is a $d$-dimensional $k$-template and   $L(\RR^d,P)$ is semialgebraically immersible into $H$, then  $L(\RR^d,P)$ is embeddable into $H$.}

\bigskip

The following theorem is the principal result of \cite{alg}.

\bigskip

{\sc Theorem} 2.4: (\cite[Theorem~2.2]{alg}) {\em Suppose that  $H$ is an algebraic $k$-hypergraph   and $\kappa$ 
is an infinite cardinal. 
 The following are equivalent$:$

\begin{itemize}

\item[$(1)$] $\chi(H) \leq \kappa;$

\item[$(2)$] if $P$ is a  $d$-dimensional $k$-template and   $L(\RR^d,P)$ is embeddable into $H$, then  $\chi(L(\RR^d,P)) \leq \kappa;$

\item[$(3)$] if  $P$ is a  $(k-1)$-dimensional $k$-template and $L(\RR^{k-1},P)$ is  semialgebraically immersible in $H$, then  $\chi(L(\RR^{k-1},P)) \leq \kappa.$

\end{itemize}}

\bigskip
  
  It should be noted that  the instance of the previous theorem when $\kappa = \aleph_0$ 
 had already appeared in \cite{avoid}. We will see in Corollary~4.3 that ``semialgebraically'' can be omitted in $(3)$ and also in Lemma~2.3.
 
   \bigskip

                   
                    {\bf \S3.\@ Compactness/Decidability.}   The  main result of this section, Theorem~3.1,  is a sort of effective compactness theorem. Various consequences of Theorem~3.1 will be presented in the  next section. 
    
     We begin with a way of constructing some new $k$-templates from an old one.       Suppose that $m \leq d < \omega$ and $\pi : d \into m$ is a surjection. If 
   $x = \langle x_0,x_1, \ldots,x_{d-1} \rangle$ is a $d$-tuple, then we define   the $\pi$-{\bf collapse} of $x$ to be the $m$-tuple $y = \langle y_0,y_1, \ldots, y_{m-1} \rangle$ where, if $j < m$, then $y_j$ 
   is the $\pi^{-1}(j)$-tuple such that if $\pi(i) = j$, then $y_{j,i}  = x_i$.  If $P$ is a $d$-dimensional 
   $k$-template, then $P^\pi$, the   $\pi$-{\bf collapse} of $P$,  is the set of the $\pi$-collapses of elements of $P$. Clearly,  $P^\pi$ is an $m$-dimensional $k$-template.  
   If $I$ is a distinguisher for $P$, then $\{\pi(j) : j \in I\}$ is a distinguisher for $P^{\pi}$, so   $e(P^\pi) \leq e(P)$.

    \bigskip

    We will say that a  polynomial  $p(x_0,x_1, \ldots, x_{k-1},y)$ is  $((k,n)+\ell)$-ary if 
    each $x_i$ is an $n$-tuple of variables and $y$ is an $\ell$-tuple of variables. 
    A formula $\varphi(z,x,y)$ is $(d+n+\ell)$-ary if $z$ is a $d$-tuple, $x$ is an $n$-tuple and $y$ is an $\ell$-tuple.
    
    \bigskip

 {\sc Theorem 3.1:}  {\em Suppose that  
 $p(x_0,x_1, \ldots, x_{k-1},y)$ is a $((k,n)+\ell)$-ary polynomial over $\QQ$  and $P$ is a $d$-dimensional $k$-template. Then there 
 are $M < \omega$ and, for each surjection $\pi : d \into m$,  an  $(m+n+ \ell)$-ary ${\mathcal L}_{OF}$-formula $\varphi_\pi(z,x,y)$ such that whenever $c \in \RR^\ell$, $H_c$ is the zero $k$-hypergraph of 
 $p(x_0,x_1, \ldots,x_{k-1},c)$ and   $L(M^d,P)$ is embeddable into $H_c$, then there  is $\pi$ such that 
 $\varphi_\pi(z,x,c)$ defines an immersion 
of $L(\RR^m,P^\pi)$ into $H_c$. 
}

\bigskip

{\it Proof}.  By replacing $p(x_0,x_1, \ldots,x_{k-1},y)$ with
 $$
 \prod_\sigma p(x_{\sigma(0)}, x_{\sigma(1)},\ldots, x_{\sigma(k-1)},y) \times \prod_{i<j<d}(x_i - x_j),
 $$
(where $\sigma$ ranges over all permutations of $k$) if needed, we can assume 
 that for each $c \in \RR^\ell$,  $p(x_0,x_1, \ldots,x_{k-1},c)$ is symmetric and reflexive.

 Let $\Pi$ be the set of all $\pi$ such that for some $m \leq d$, $\pi : d \into m$ is a surjection. 
 If $\pi \in \Pi$ and $m$ is the range of $\pi$, then $m = \pi[d]$.
 When we consider a formula $\varphi_\pi(z,x,y)$ then it is to be understood that it is a  
 $(\pi[d]+n+\ell)$-formula.

Suppose, for a contradiction, that there are no such $M$ and $\langle\varphi_\pi(z,x,y) : \pi \in \Pi \rangle$ as in the theorem.   Thus, we have the following:

\begin{quote} 
\hspace{-26pt} ($*$) For every $M < \omega$ and $\langle\varphi_\pi(z,x,y) : \pi \in \Pi \rangle$, there is $c \in \RR^\ell$ such that $L(M^d,P)$ is embeddable into $H_c$ and for 
 no $\pi \in \Pi$ does 
 $\varphi_\pi(z,x,c)$  define an immersion 
of $L(\RR^{\pi[d]},P^\pi)$ into $H_c$. 
\end{quote}

\noindent The statement $(*)$ implies the following stronger one: 

\begin{quote}
\hspace{-31pt} ($**$) For every $M < \omega$ and finitely many 
$\langle \varphi_{0,\pi}(z,x,y) : \pi \in \Pi \rangle, 
\langle \varphi_{1,\pi}(z,x,y) : \pi \in \Pi \rangle,  \ldots,  \langle \varphi_{N,\pi}(z,x,y) : \pi \in \Pi \rangle$,  there is $c \in \RR^\ell$ such that 
$L(M^d,P)$ is embeddable into $H_c$ and for no $i \leq N$ and $\pi \in \Pi$ does $\varphi_{i,\pi}(z,x,c)$   define an immersion 
of $L(\RR^{\pi[d]},P^\pi)$ into $H_c$. 
\end{quote}

To prove $(**)$, suppose that $M < \omega$ and 
$\langle \varphi_{i,\pi}(z,x,y): i \leq N,\pi \in \Pi  \rangle$ constitute a counterexample. 
Thus, for each $c \in \RR^\ell$, either $L(M^d,P)$ is not embeddable in $H_c$ or there are $i \leq N$ and $\pi \in \Pi$ such that $\varphi_{i,\pi}(z,x,c)$ defines an immersion of 
$L(\RR^{\pi[d]},P^{\pi})$ into $H_c$. For each $\pi \in \Pi$ and $i \leq N$, let $C_{i,\pi} $ be the set of all $c \in \RR^\ell$ such that 
 $H_c$ does not embed $L(M^d,P)$ or $\varphi_{i,\pi}(z,x,c)$   defines an immersion 
of $L(\RR^\pi[d],P^\pi)$ into $H_c$. Each $C_{i,\pi}$ is semialgebraic and $\RR^\ell = \bigcup_{i \leq N} \bigcup_{\pi \in \Pi}C_{i,\pi}$. 
 For each $\pi \in \Pi$,  let $\varphi_\pi(z,x,y)$ be the formula
$$
 \bigvee_{i \leq N}\varphi_{i,\pi}(z,x,y) \wedge y \in C_{i,\pi} \backslash \bigcup_{j < i} C_{j,\pi}.
 $$
 We then have that  for every $c \in \RR^\ell$, if $L(M^d,P)$ is embeddable into $H_c$, then there is $\pi \in \Pi$ such that $\varphi_\pi(z,x,c)$  defines an immersion 
of $L(\RR^{\pi[d]},P^\pi)$ into $H_c$. This contradicts $(*)$ and, thereby, proves $(**)$.

\smallskip

 Let $\lambda > \beth_{\omega+\omega}$ be a  cardinal such that $\lambda^{<\lambda} = \lambda$.
(There is no guarantee that such a cardinal exists; if there is none, then work inside an 
appropriate inner model of the universe of sets that has such a cardinal.) Utilizing  these properties
of $\lambda$, we let 
$ R$ be a saturated elementary extension of {$ \RR$} 
such that $|R| = \lambda$.

Since ${ R}$ is saturated,  $\lambda > \beth_{\omega+\omega}$ and $(**)$ holds, we can get $c  \in R^\ell$ such that (letting $H$ be the zero $k$-hypergraph of 
$p(x_0,x_1, \ldots, x_{k-1},c)$ in $R$) such that:   
 
 \begin{itemize}

\item[$\circledast_1$] $L((\beth_{\omega+\omega})^d,P)$ is embeddable into $H$.

\item[$\circledast_2$] For each $\pi \in \Pi$, there is  no formula $\varphi_\pi(z,x,y)$ such that 
$\varphi_\pi(z,x,c)$ defines in $R$ and immersion of $L(R^{\pi[d]},P^\pi)$  into $H$. 

\end{itemize}

 We let $c$ and  $H$ be fixed for the rest of this proof. Rewording $\circledast_2$, we have that  
 for each $\pi \in \Pi$, there is no  $\{c_0,c_1, \ldots,c_{\ell-1}\}$-definable immersion of $L(R^{\pi[d]},P^\pi)$ into $H$. Since $H$ is $\{c_0,c_1, \ldots ,c_{\ell-1}\}$-definable, if there were an $R$-semialgebraic 
immersion $L(R^{\pi[d]},P^\pi)$  into $H$, then there would be one that is $\{c_0,c_1, \ldots,c_{\ell-1}\}$-definable in $R$ since ${\mathsf {RCF}}$ has definable Skolem functions 
\cite{vdd}. Thus, we can strengthen $\circledast_2$ to:
\begin{itemize}

\item[$\circledast_3$]  For each $\pi \in \Pi$, there is no $R$-semialgebraic immersion of $L(R^{\pi[d]},P^\pi)$  into $H$.

\end{itemize}

  Let $T$ be a transcendence basis for $R$.  Because $R$ is saturated,
  we can require that $T$ be dense;  that is, whenever $a,b \in R$ and $a < b$, 
  then $(a,b) \cap T \neq \varnothing$. 
   To get such a $T$, 
  first by transfinite recursion of length $\lambda$, construct a 
   a dense $T_0 \subseteq R$ that is algebraically independent, and  then  extend it to a transcendence basis $T$.  
  Because $R$ is saturated, it follows that $|(a,b) \cap T| = \lambda$ whenever $a < b \in R$.

   For the next definitions, suppose that  $F \subseteq T$ is finite.

     For each $r < \omega$ and $a = \langle a_0,a_1, \ldots, a_{r-1} \rangle  \in R^r$, define $\supp_F(a)$, 
     the $F$-{\bf support} of $a$,   to be the smallest subset $S \subseteq T$ such 
     that $\{a_0,a_1, \ldots, a_{r-1}\}$ is an $(S \cup F)$-definable subset of $R$. Equivalently, 
     $\supp_F(a)$ is the smallest $S \subseteq T$ such that each $a_i$ is algebraic over $S  \cup F$. 
  Easily, $\supp_F(a)$ is a well defined, finite 
  subset of $T \backslash F$. 
  
  Let $F_0 = \supp_\varnothing(c)$. 
    
    We say that $f$ is an $F$-{\bf determining function} for  $a \in R^r$, where  
    $\supp_F(a) = \{t_0,t_1, \ldots,t_{m-1}\}_< \subseteq T^m$, if the following hold:         
    
    \begin{itemize} 
    
    \item $f : U \into R^n$  is $F$-definable in $R$;
    
    \item    $U \subseteq R^{(m)}$ and there is an $R$-Nash homeomorphism from $R^m$ onto $U$;

    \item $f$ is one-to-one on each coordinate;

    \item $\langle t_0,t_1, \ldots,t_{m-1}\rangle \in U$ and $f(t_0,t_1, \ldots, t_{m-1}) = a$.
     
    \end{itemize}
    
    \bigskip

    {\sc Lemma 3.1.0}: {\em For each $r < \omega$ and $a \in R^r$, there is an $F$-determining function. In fact, there is one that is Nash.}
    
    \bigskip
    
    {\it Proof.} (Sketch) First, suppose that $r=1$ so that $a \in R$. Let $\supp_F(a) = \{t_0,t_1, \ldots, t_{m-1}\}_<$. Let $p(x_0,x_1, \ldots, x_{m-1},y)$ be in irreducible polynomial with coefficients in the real closed field generated by $F$ such that $p(t_0,t_1, \ldots,t_{m-1}, a) = 0$. Suppose that $a$ is the {\mbox {$k$-th }} largest root. Let $g(x_0,x_1, \ldots, x_{m-1}) = y$ be such that $y$ is the $k$-th largest root of $p(x_0,x_1, \ldots, x_{m-1},y)$. Then by use of a cell stratification \cite[Chap.\@ 5.4]{arag}, there is $U \subseteq R^m$ so that $f = g \harp U$ is as required.
    
    Now let $a \in R^r$. Let $\supp_F(a) = \{t_0,t_1, \ldots, t_{m-1}\}_<$. For each $i < r$, let $J_i \subseteq m$ be such that $\supp_F(a_i) = \{t_j : j \in J_i\}$. Let $f_i : U_i \into R$ be an $F$-determining function for $a_i$. (We arrange that $U_i \subseteq R^{J_i}$.) For each $i < m$, let 
    $p_i : R^r \into R^{J_i}$ be the projection map. Let $U \subseteq R^{(m)}$   be the largest set such 
    $p_i[U] \subseteq U_i$ for $i < r$. Let $f : U \into R^r$ be such that $f(x)_i = f_ip_i(x)$. \qed
    
    \bigskip

   Since $L((\beth_{\omega+\omega})^d,P)$ is embeddable in $H$, there certainly is an  embedding of $L(\omega^d,P)$ into $H$. 
      In Lemma~3.1.2, we  will obtain  an embedding of $L(\omega^d,P)$ into $H$ having some additional properties.

    \bigskip
    
    {\sc Lemma 3.1.1}: {\em There are  $f$ and  an embedding $h$ of  $L\big((\beth_\omega)^d,P\big)$ into $H$  such that 
     whenever $\alpha \in (\beth_\omega)^d$, then $f$ is an $F_0$-determining function for $h(\alpha)$.}
     
      \bigskip 
    
    {\it Proof}. Let $\theta$ be an embedding of $L\big((\beth_{\omega + \omega})^d,P\big)$ into $H$. Consider a  function $\varphi$ on $(\beth_{\omega + \omega})^d$ such that 
    if $\alpha \in (\beth_{\omega+\omega})^d$, then $\varphi(\alpha)$ is 
    a $F_0$-determining function for $\theta(\alpha)$.  The range of $\varphi$ has cardinality at most ${\aleph_0}$   since each $f$ in the range is   $F_0$-definable in ${ R}$. Thus, by Corollary~1.2,  there are  $f$ and 
    $Y_0,Y_1, \ldots, Y_{d-1} \subseteq \beth_{\omega+\omega}$  such that $|Y_0| = |Y_1| = \cdots = |Y_{d-1}| = \beth_\omega$ and $\varphi$ is constantly $f$ on $Y = Y_0 \times Y_1 \times \cdots \times Y_{d-1}$. Thus, $\theta \harp Y$ is an embedding of $L(Y,P)$ into $H$. Since $L((\beth_\omega)^d,P) \cong L(Y,P)$,   this implies the existence of the required~$h$. \qed

    \bigskip
  
Suppose that $F_0 \subseteq F \subseteq T$ and that $F$ is finite, and let $h : X^d \into R^n$ be an embedding of $L(X^d,P)$ into $H$.   We define 
the function $\sigma_{h,F}$ on $X^d$ so that if $\alpha \in X^d$ and 
$\supp_F(h(\alpha)) = \{t_0,t_1, \ldots, t_{m-1} \}_<$,       then $\sigma_{h,F}(\alpha) = \langle t_0,t_1, \ldots, t_{m-1} \rangle$. Thus, 
    if $h$ and $f$ are  such that for each $\alpha \in X^d$, $f$ is an $F$-determining function of $h(\alpha)$ and 
    $\dom(f) \subseteq R^m$, then $\sigma_{h,F}(\alpha) \in (T \backslash F)^{(m)}$ and 
    $h(\alpha) = f\big(\sigma_{h,F}(\alpha)\big)$ for every $\alpha \in X^d$.

    \bigskip
    
    {\sc Lemma 3.1.2}: {\em Let $F = F_0$. There are $m < \omega$ and  an $m$-ary  function $f$ such that for each $\gamma < \beth_{\omega}$, there are an embedding $h $ of $L(\gamma^d,P)$  
    into $H$ and a function $D : {\mathcal P}(m) \into {\mathcal P}(d)$ such that$:$
    
      \begin{itemize}

    \item[Z$_1(\gamma)\hspace{-3pt}:$]  whenever $\alpha \in \gamma^d$, then $f$ is an $F$-determining function for $h(\alpha);$
    
      \item [Z$_2(\gamma)\hspace{-3pt}:$]  whenever 
     $I \subseteq m$ and $\alpha, \beta \in \gamma^d$, then
     $$
     \alpha \sim_{D(I)} \beta \Longleftrightarrow \sigma_{h,F}(\alpha) \sim_I \sigma_{h,F}(\beta);
     $$
    
    \item[Z$_3(\gamma)\hspace{-3pt}:$]  whenever $\alpha, \beta \in \gamma^d$ and $i < j < m$, 
    then $\sigma_{h,F}(\alpha)_i \neq \sigma_{h,F}(\beta)_j$. 
 
       \end{itemize}}
       
       \bigskip
       
       {\it Proof}. Let $f$ be  as in Lemma~3.1.1, and  let $m < \omega$ be such that $f$ is $m$-ary. Thus, we already know (Lemma~3.1.1) that for every $\gamma < \beth_\omega$, there is an 
       embedding $h$ of  $L(\gamma^d,P)$ into $H$ satisfying Z$_1(\gamma)$. 
       
       To take care of Z$_2(\gamma)$, let $I_0,I_1, \ldots, I_{2^m-1}$ be all the subsets of $m$. 
       For each $i \leq 2^m$, let $({\mathsf S}_i)$ be the statement: 
       
       \begin{quote}
       
       $({\mathsf S}_i)$: For every $\gamma < \beth_\omega$, there is an 
       embedding $h$ of  $L(\gamma^d,P)$ into $H$ satisfying Z$_1(\gamma)$ and such that 
       there are $J_0,J_1, \ldots, J_{i-1} \subseteq d$  such that whenever $j < i$  
       and $\alpha,\beta \in \gamma^d$, then 
       $\alpha \sim_{J_j} \beta \Longleftrightarrow \sigma_h(\alpha) \sim_{I_j} \sigma_h(\beta)$.
       
       \end{quote} \noindent
       We  prove $({\mathsf S}_i)$ be induction on $i \leq 2^m$. The basis step $({\mathsf S}_0)$ is essentially vacuously true.
       Given that $({\mathsf S}_i)$ is true, we can prove $({\mathsf S}_{i+1})$ by  applications of Theorem~1.3.
       Having that $({\mathsf S}_{2^m})$ is true, we just let $D(I_j) = J_j$ for $j < 2^m$
       
       Finally, to get Z$_3(\gamma)$, just apply  Corollary~1.5. \qed
       
       \bigskip

    For the remainder of the proof of Theorem~3.1,  we fix $m < \omega$, a finite $F \subseteq T$ such that $F \supseteq F_0$ and an $F$-determining function $f$ for which  there are an embedding $h$ of  $L(\omega^d,P)$  
    into $H$ and a function $D : {\mathcal P}(m) \into {\mathcal P}(d)$ such that$:$

        \begin{itemize}

    \item[Z$_1(\omega)\hspace{-3pt}:$]  whenever $\alpha \in \omega^d$, then $f$ is an $F$-determining function for $h(\alpha);$
    
      \item [Z$_2(\omega)\hspace{-3pt}:$]  whenever 
     $I \subseteq m$ and $\alpha, \beta \in \omega^d$, then
     $$
     \alpha \sim_{D(I)} \beta \Longleftrightarrow \sigma_{h,F}(\alpha) \sim_I \sigma_{h,F}(\beta);
     $$
    
    \item[Z$_3(\omega)\hspace{-3pt}:$]  whenever $\alpha, \beta \in \omega^d$ and $i < j < m$, then 
    $\sigma_{h,F}(\alpha)_i \neq \sigma_{h,F}(\beta)_j$. 
    
    \end{itemize}
    \smallskip
     Furthermore, we require that $m$, $F$ and $f$ are chosen so as to minimize~$m$. 
     Lemma~3.1.2 assures that this is possible. We also let $U = \dom(f)$. From now on we drop the ``$F$''  and write ``determining'', ``supp'' and ``$\sigma_h$" instead of 
``$F$-determining'', ``$\supp_F$'' and ``$\sigma_{h,F}$", respectively.

    \bigskip   
    
 We will say that 
    $h$ is a $D$-{\bf normal embedding} if 
    $h$ is an embedding of $L(\omega^d,P)$ into $H$ and $D : {\mathcal P}(m) \into {\mathcal P}(d)$ is such that Z$_1(\omega)$--Z$_3(\omega)$  hold.  Let ${\mathcal D}$ be the set of all 
    $D : {\mathcal P}(m) \into {\mathcal P}(d)$ for which  there is a 
    $D$-normal  embedding.  Lemma~3.1.2 implies that ${\mathcal D} \neq \varnothing$. 
    
       The next lemma states some  properties of those $D \in {\mathcal D}$.
    
   \bigskip
   
   {\sc Lemma 3.1.3}: {\em  Suppose that $D \in {\mathcal D}$ and $I,J \subseteq m$. Then$:$
   
   \begin{itemize}
   
   \item[(W$_1$)] $D(I) = \varnothing$ iff  $I = \varnothing;$
   
    \item[(W$_2$)] $D(I) = d$ if $I = m;$

   \item[(W$_3$)] $D(I \cup J) =  D(I) \cup D(J)$.

   \end{itemize}}
   
   \bigskip

   {\it Proof}.    Suppose that $D \in {\mathcal D}$ and $I,J \subseteq m$. 
   
   \smallskip
   
   (W$_2$):  Since $h$ is one-to-one, then $\sigma_h$ is also one-to-one implying that $D(m) = d$.

\smallskip
   
   (W$_3$): We have that 
      \begin{eqnarray*}
   \alpha \sim_{D(I \cup J)} \beta & \Longleftrightarrow  & \sigma_h(\alpha) \sim_{I \cup J} \sigma_h(\beta)
   \\
   & \Longleftrightarrow & \sigma_h(\alpha) \sim_I \sigma_h(\beta) \mbox{ and } \sigma_h(\alpha) \sim_J \sigma_h(\beta) \\
   & \Longleftrightarrow & \alpha \sim_{D(I)} \beta \mbox{ and } \alpha \sim_{D(J)} \beta \\
    & \Longleftrightarrow & \alpha \sim_{D(I) \cup D(J)} \beta.
     \end{eqnarray*}

     \smallskip

  (W$_1$): If $I = \varnothing$, then it is obvious that $D(I) = \varnothing$. 
  To prove the converse implication, assume, for a contradiction, that $I \neq \varnothing$ and that 
  $D(I) = \varnothing$. By (W$_3$), we can assume that $I = \{i\}$, where $i < m$.    
  
      Let $h$ be a $D$-normal embedding. 
   For any $\alpha, \beta \in \omega^d$, $\sigma_h(\alpha)_i = \sigma_h(\beta)_i$. Thus, there is $t \in 
   T \backslash F$ such that $\sigma_h(\alpha)_i = t$ for each $\alpha \in \omega^d$. Let $U_0 = \{\langle t_0,t_1, \ldots,t_{i-1},t_{i+1}, \ldots, t_{m-1} \rangle \in T^{m-1} : \langle t_0,t_1, \ldots,t_{i-1},t,t_{i+1}, \ldots, t_{m-1} \rangle \in U\}$. Let $f_0 : U_0 \into R^n$ be such that 
   $$
   f_0( t_0,t_1, \ldots,t_{i-1},t_{i+1}, \ldots, t_{m-1}) = f(t_0,t_1, \ldots,t_{i-1},t,t_{i+1}, \ldots, t_{m-1}).
   $$
   It is clear that $f_0$ is an $(F \cup\{t\})$-determining function of $h(\alpha)$ for each $\alpha \in \omega^d$.  Clearly, $f_0$ is $D_0$-normal, where 
   $D_0(J) = D((J \cap i) \cup \{j+1: i \leq j \in D\})$ for $J \subseteq m-1$, 
   and also 
    $\sigma_{h,F \cup \{t\}}(\alpha)_i \neq \sigma_{h,F \cup \{t\}}(\beta)_j$ whenever $\alpha, \beta \in \omega^d$ and $i < j < m-1$.

   This contradicts the minimality of $m$.  \qed

   \bigskip
   
    A consequence of (W$_1$) and (W$_3$)  is that each $D \in {\mathcal D}$  is completely determined by its value on singletons. Also, (W$_3$) has a consequence that $D(I \cap J) \subseteq D(I) \cap D(J)$.

    \bigskip

   {\sc Lemma 3.1.4}: {\em Suppose that  $D \in {\mathcal D}$ and that 
    $B \subseteq U$ is an open $m$-cube. Then there is  a $D$-normal embedding $h$ such that 
   $\sigma_h : \omega^d \into B$.}
   
   \bigskip
   
   {\it Proof}.  Suppose that $B = (a_0,b_0) \times (a_1,b_1) \times \cdots \times (a_{m-1},b_{m-1}) \subseteq U$ is an  open $m$-cube. Since $U \subseteq R^{(m)}$, then   $a_0 < b_0 \leq a_1 < b_1 \leq  \cdots < b_{m-2} \leq a_{m-1} < b_{m-1}$. 
   Let $\theta$ be a $D$-normal embedding of $L(\omega^d,P)$. 
   Let $S_0,S_1, \ldots,S_{m-1} \subseteq T \backslash F$ be the smallest sets such that 
   $\sigma_\theta : \omega^d \into S_0 \times S_1 \times \cdots \times S_{m-1}$. Each $S_i$ is countable.
   Since $T$ is dense, there are one-to-one functions $g_i : S_i \into (a_i,b_i) \cap (T \backslash F)$ for  $i < m$. 
   Let $g : S_0 \times S_1 \times \cdots \times S_{m-1} \into B$ be such that 
   $g(\langle t_0,t_1, \ldots, t_{m-1} \rangle) = \langle g_0(t_0),g_1(t_1), \ldots, g_{m-1}(t_{m-1}) \rangle$.
  We can now define $h : \omega^d \into R^n$  so that 
  $h(\alpha) = fg\sigma_\theta(\alpha)$. We claim that $h$ is as required. Since it is the composition of one-to-one functions,  $h$ is one-to-one.
  
  \smallskip
  
  {\em $h$ is an embedding of $L(\omega^d,P)$ into $H$.}  Let $\{\alpha_0,\alpha_1, \ldots, \alpha_{k-1}\}$ be an edge of $L(\omega^d,P)$. Then $\{\theta(\alpha_0), \theta(\alpha_1), \ldots, \theta(\alpha_{k-1})\}$ is an edge of $H$. For $i < k$, let $\sigma_\theta(\alpha_i) = t_i \in (T \backslash F)^m$. Thus,  $p(f(t_0),f(t_1), \ldots, f(t_{k-1})) = 0$, so  $p(f(g(t_0)),f(g(t_1)), \ldots, f(g(t_{k-1}))) = 0$ by   Lemma~1.6. Therefore, 
  $\{h(\alpha_0),h(\alpha_1), \ldots, h(\alpha_{k-1})\}$ is an edge of $H$, so $h$ is an embedding.
  
  \smallskip
  
  {\em $h$ is $D$-normal.} Clearly, whenever $I \subseteq m$ and $\alpha, \beta \in \omega^d$, 
  then $\sigma_\theta(\alpha) \sim_I \sigma_\theta(\beta)$ iff  $\sigma_\theta(g(\alpha)) \sim_I \sigma_\theta(g(\beta))$. But $g(\alpha) = \sigma_h(\alpha)$, so that 
  $\sigma_\theta(\alpha) \sim_I \sigma_\theta(\beta)$ iff  $\sigma_h(\alpha) \sim_I \sigma_h(\beta)$. Hence, since $\theta$ is $D$-normal, then so is $h$. 
 \qed
 
   \bigskip

  We will say that $D \in {\mathcal D}$ is  {\bf minimal} if  whenever $D' \in {\mathcal D}$ is 
   such that $D'(I) \subseteq D(I)$ for every $I \subseteq m$, then   $D' = D$. 
   Since ${\mathcal D}$ is finite, there is at least one minimal  $D \in {\mathcal D}$.      
          
   \bigskip
   
   {\sc Lemma~3.1.5}: {\em Suppose that  $D \in {\mathcal D}$ is minimal. If $i < j < m$, then $D(\{i\}) \cap D(\{j\}) = \varnothing$.}
   
   \bigskip
   
   {\it Proof}. With the aim of contradicting that $D$ is minimal, we assume that $i < j < m$ and that $D(\{i\}) \cap D(\{j\}) \neq \varnothing$. 
   It must be that $m \geq 2$. Without loss of generality, let $0 \in D(\{i\}) \cap D(\{j\})$. 
   
  Let $B = B_0 \times B_1 \times \cdots \times B_{m-1}  \subseteq U$ be an open $m$-cube. By Lemma~3.1.4,    let $h$ be a  $D$-normal embedding such that $\sigma_h : \omega^d \into B$.  
   We wish to modify $\sigma_h$ so as   to get $\tau : \omega^d \into B \cap (T \backslash F)^m$.
   
   Consider $\alpha \in \omega^d$. To define $\tau(\alpha)$, we will define $\tau(\alpha)_\ell \in B_\ell\cap (T \backslash F)$ for 
   each $\ell < m$. First, if $i \neq \ell < m$, then let $\tau(\alpha)_\ell = \sigma_h(\alpha)_\ell$. 
   Next, let $\tau(\alpha)_i = \sigma_h(\langle 0,\alpha_1,\alpha_2, \ldots, \alpha_{d-1} \rangle)_i \in 
   B_i \cap (T \backslash F)$. Thus, $\tau(\alpha)_\ell \in B_\ell \cap (T \backslash F)$ for each $\ell < m$ so that 
   $\tau(\alpha) \in B \cap (T \backslash F)^m$.  
   
   We claim that $\tau$ is one-to-one. To see this, consider distinct $\alpha, \beta \in \omega^d$ intending to show that $\tau(\alpha) \neq \tau(\beta)$.
   If  $i \neq \ell < m$ and $\sigma_h(\alpha)_\ell \neq   \sigma_h(\beta)_\ell$, then 
   $\tau(\alpha)_\ell = \sigma_h(\alpha)_\ell \neq   \sigma_h(\beta)_\ell = \tau(\beta)_\ell$.  
    Thus, we can assume that 
   $\sigma_h(\alpha)_\ell = \sigma_h(\beta)_\ell$ whenever $i \neq \ell < m$. Since $i \neq j$, then 
   $\sigma_h(\alpha)_j = \sigma_h(\beta)_j$, so that $\alpha_r = \beta_r$ for all $r \in D(\{j\})$ and, in particular,  $\alpha_0 = \beta_0$. Since $\sigma_h$ 
   is one-to-one, we have that $\sigma_h(\alpha)_i \neq \sigma_h(\beta)_i$. 
   We will conclude that $\tau(\alpha)_i \neq \tau(\beta)_i$. For, suppose to the contrary that 
   $\tau(\alpha)_i = \tau(\beta)_i$. Then
   $\sigma_h(\langle 0,\alpha_1,\alpha_2, \ldots, \alpha_{d-1} \rangle)_i = \sigma_h(\langle 0,\beta_1, \beta_2, \ldots, \beta_{d-1} \rangle)_i$. But then $\alpha_r = \beta_r$ whenever $0 < r \in D(\{i\})$. But, since $\alpha_0 = \beta_0$, we get that $\sigma_h(\alpha)_i = \sigma_h(\beta)_i$, which is a contradiction, proving that $\tau$ is one-to-one.
   
  Let $g = f \tau$, so that $g : \omega^d \into R^n$.  Since $\tau$ is one-to-one, then also $g$ is one-to-one.
   
   We claim that $g$ is an embedding of $L(\omega^d,P)$ into $H$.  Suppose that $\{\alpha_0,\alpha_1, \ldots, \alpha_{k-1}\}$ is an edge of $L(\omega^d,P)$, intending to show that 
   $\{g(\alpha_0), g(\alpha_1), \ldots,g(\alpha_{k-1})\}$ is an edge of $H$ or, equivalently,
   that 
   $$
   p\big(g(\alpha_0), g(\alpha_1), \ldots, g(\alpha_{k-1})\big) = 0.
   $$
   For each $r < k$, let $\sigma_h(\alpha_r) = t_r = \langle t_{r,0}, t_{r,1}, \ldots, t_{r,m-1} \rangle \in (T \backslash F)^m$. 
   Since $h$ is an embedding, we have that
   $
   p\big(f(t_0),f(t_1), \ldots, f(t_{d-1})\big) = 0.
   $
   For each $r <k$, let $\tau(\alpha_r) = t'_r = \langle t'_{r,0}, t'_{r,1}, \ldots, t'_{r,m-1} \rangle \in (T \backslash F)^m$.
   If $i \neq \ell < m$, then $t'_{r,\ell} = t_{r,\ell}$.  If $r < s < m$ and $t_{r,i} = t_{s,i}$, then $t'_{r,i} = t'_{s,i}$. 
   By Z$_3(\omega)$, if $r,s < k$ and $j < \ell < m$, then $t_{rj} \neq t_{s\ell}$ and $t'_{rj} \neq t'_{s\ell}$. 
   Since $p$ and $f$ are semialgebraic and analytic, we have that  
   $p\big(f(t'_0),f(t'_1), \ldots, f(t'_{d-1})\big) = 0.$
   Thus, $f\tau$ is an embedding $L(\omega^d,P)$ into $H^R$. 
   
   One checks that $\tau = \sigma_g$. 
   
   Let $D' : {\mathcal P}(m) \into {\mathcal P}(d)$ be such that if $I \subseteq m$, then 
   \begin{displaymath}
   D'(I) = \left\{  \begin{array}{ll}
   D(I)  & \textrm{if }\  i \not\in I \\
   D(I) \backslash \{0\} & \textrm{if }\ i \in I.
   \end{array}  \right.
   \end{displaymath}
   One easily checks that whenever $I \subseteq m$ and $\alpha, \beta \in \omega^d$, then 
   $\alpha \sim_{D'(I)} \beta$ iff $\tau(\alpha) \sim_I \tau(\beta)$, thereby contradicting the minimality of $D$.   \qed
   
   \bigskip

   A consequence of Lemma~3.1.5  is that if $D \in {\mathcal D}$ is minimal and $I,J \subseteq m$, then the aforementioned consequence of (W$_3$)  can be improved to: {\em     $D(I \cap J) = D(I) \cap D(J)$.}

     \smallskip
     
       We now fix a minimal  $D \in {\mathcal D}$.   Let $\pi : d \into m$ be the unique (by Lemmas~3.1.3 and~3.1.5)  function such that if $ j < d$, then  $j \in D(\{\pi(j)\})$.

 \bigskip
 
 {\sc Lemma 3.1.6}: {\em There is an $R$-semialgebraic immersion of $L(R^m,P^\pi)$ into $H$.}
 
 \bigskip
 
 {\it Proof}. Let $B = B_0 \times B_1 \times \cdots \times B_{m-1} \subseteq U$ be an open $m$-cube. 
  It suffices to show that there is an $R$-semialgebraic immersion of $L(B,P\pi)$ into $H$. We will prove, in fact, that $f \harp B$ is such an immersion. Since $f$ is  a determining function, it is an $R$-Nash function that is  one-to-one in each coordinate. Thus, we will complete the proof upon proving:
 {\em If $\{x_0,x_1, \ldots,x_{k-1}\}$ is an edge of $L(B,P^\pi)$ and $f(x_0),f(x_1), \ldots, f(x_{k-1})$ are pairwise distinct, then $\{f(x_0),f(x_1), \ldots, f(x_{k-1})\}$ is an edge of $H$.} 
 Or, equivalently:  {\em If $\{x_0,x_1, \ldots,x_{k-1}\}$ is an edge of $L(B,P^\pi)$, then $p(f(x_0),f(x_1), \ldots, f(x_{k-1}))=0$.} But then it is enough to show that there is some $m$-dimensional $k$-template 
 $\{t_0,t_1, \ldots,t_{k-1}\} \subseteq T^m \cap B$ such that $P^\pi$ is a  homomorphic image of 
 $\{t_0,t_1, \ldots,t_{k-1}\}$ and 
$p(f(t_0),$ $f(t_1), \ldots, f(t_{k-1}))=0$.
  
 By Lemma~3.1.4, let $h$ be a $D$-normal embedding such that $\sigma_h : \omega^d \into B$. 
 Let $\{\alpha_0,\alpha_1, \ldots,\alpha_{k-1}\}  \subseteq \omega^d$ be a $d$-dimensional 
 $k$-template isomorphic to $P$. Without loss, assume that $P = \{\alpha_0, \alpha_1, \ldots, \alpha_{k-1}\}$. Then, $p(h(\alpha_0),h(\alpha_1), \ldots, h(\alpha_{k-1})) = 0$. 
 Let $t_i = \sigma_h(\alpha_i)$. Thus,  $\{t_0,t_1,$ $ \ldots,t_{k-1}\} \subseteq T^m \cap B$ and 
$p(f(t_0),f(t_1), \ldots, f(t_{k-1}))=0$. It remains to show that $P^\pi$ is a homomorphic image of 
$\{t_0,t_1, \ldots,t_{k-1}\}$. 

For $r < k$, let $x_r$ be the $\pi$-collapse of $\alpha_r$. 
Consider $i < m$ and $r,s < k$, and suppose that $t_{r,i} = t_{s,i}$. Then $\alpha_r \sim_{D(\{i\})}  \alpha_s$, so that $\alpha_{r,j} = \alpha_{s,j}$ whenever $\pi(j) = i$. Thus, $x_{r,i} = x_{s,i}$.
\qed

\bigskip

Lemma~3.1.6 contradicts $\circledast_3$, thereby completing the  proof of Theorem~3.1. \qed     \bigskip


    {\bf \S4.\@ Algebraic hypergraphs.}  Corollary~4.3, our first corollary to Theorem~3.1, is an improved compactness-type theorem for algebraic hypergraphs.  Before proving it, we need two lemmas, the first of which is a simple combinatorial lemma.
    
    \bigskip
    
    {\sc Lemma 4.1}: {\em Suppose that $d < \omega$ and $X$ is an infinite set. 
    If $f : X^d \into C$ is one-to-one in each coordinate, then, for each $M < \omega$, there are 
    $Y_0,Y_1, \ldots, Y_{d-1} \in [X]^M$ such that  $f$ is one-to-one on $Y_0 \times Y_1 \times \cdots \times Y_{d-1}$.}
    
    \bigskip
    
    {\it Proof.} Let $E$ be the finite set of all equivalence relations on $\{0,1\}^d$. 
    Consider some $M < \omega$,  and assume that $M \geq 3$. 
    Without loss of generality, let  $X = \omega$.   Define a function $p : ([\omega]^2)^d  \into E$ as follows. 
  
 Suppose that $A = \langle A_0,A_1, \ldots,A_{d-1} \rangle \in ([\omega]^2)^d$, where, for each $i < d$,  $A_i = \{a_{i,0}, a_{i,1}\}$ and $a_{i,0} < a_{i,1}$.     
 Now we define $p(A) = \theta \in E$, where $\theta$ is  such that whenever $\sigma, \tau \in \{0,1\}^d$, then 
 $$
 \langle \sigma,\tau \rangle \in \theta \Longleftrightarrow f(a_{i,\sigma_0},a_{i,\sigma_1}, \ldots, a_{i,\sigma_{d-1}}) = 
    f(a_{i,\tau_0},a_{i,\tau_1}, \ldots, a_{i,\tau_{d-1}}).
    $$
    
    By  a Polarized Ramsey's Theorem, which is an easy consequence of Ramsey's Theorem, get 
    $Y_i \in  [\omega]^M$ for each $i < d$ such that $p$ is constant on 
    $[Y_0]^2 \times [Y_1]^2 \times \cdots \times [Y_{d-1}]^2$.
    
    We  claim that $f$ is one-to-one on $Y_0 \times Y_1 \times \cdots \times Y_{d-1}$. 
    To see this, consider distinct $y,z \in Y_0 \times Y_1 \times \cdots \times Y_{d-1}$ and suppose, for a contradiction, that $f(y) = f(z)$. Let $i < d$ be such that $y_i \neq z_i$,  and assume that $y_i < z_i$.
     Since $|Y_i| = M \geq 3$, there is $c \in Y_i$ such that either 
    $y_i \neq c < z_i$ or $y_i < c \neq z$. (For definitiveness, let $c$ satisfy the former.) 
     Now let $A = \langle A_0,A_1, \ldots,A_{d-1} \rangle \in [Y_0]^2 \times [Y_1]^2 \times \cdots \times [Y_{d-1}]^2$ be such that $y,z \in A_0 \times A_1 \times \cdots \times A_{d-1}$. 
     Thus,  $A_i = \{y_i,z_i\}$. Let $B =  B_0 \times B_1 \times \cdots \times B_{d-1}$ be such that 
     $B_i = \{c,z_i\}$ and $B_j = A_j$ when $i \neq j < d$. Let $x = \langle x_0,x_1, \ldots,x_{d-1} \rangle$ be such that $x_i = c$ and $x_j = y_j$ when $i \neq j < d$. Thus, $x,z \in B$. Since $p$ is constant on  
     $[Y_0]^2 \times [Y_1]^2 \times \cdots \times [Y_{d-1}]^2$, then $f(x) = f(z)$. But then $f(x) = f(y)$,
     contradicting that $f$ is one-to-one in the $i$-th coordinate. \qed
     
     \bigskip
     
     The second lemma conerns embedding one template hypergraph into  another.
     
     \bigskip
     
     {\sc Lemma 4.2}: {\em Suppose that $m \leq d < \omega$ and $\pi : D \into m$ is a surjection. Let $P$ be a $d$-dimensional $k$-template. Then, for any infinite set $X$, $L(X^d,P)$ is embeddable into $L(X^m,P^\pi)$. $\big($Morover, $L(X^d,P)$ is isomorphic to a spanning subhypergraph of $L(X^m,P^\pi).\big)$}
     
     \bigskip
     
     {\it Proof.} For each $i < m$, let $J_i = \pi^{-1}(i)$. Since $\pi$ is surjective, each $J_i  \neq \varnothing$. Thus, $\{J_0,J_1, \ldots, J_{m-1}\}$ is a partition of $d$ into $m$ parts. Let $Y_i = X^{J_i}$. Then $|Y_i| = |X^{J_i}| = |X|$ since $X$ is infinite. Thus, it suffices to get an embedding of $L(X^d,P)$ into $L(Y,P^\pi)$, where $Y = Y_0 \times Y_1 \times \cdots \times Y_{m-1}$.
     
     Assume that $P \subseteq X^d$. Let $f_\pi : X^d \into Y$ be the $\pi$-collapse. Thus, $f_\pi[P] = P^\pi$. Clearly, $f_\pi$ is a bijection. We check that $f_\pi$ is 
     an embedding of $L(X^d,P)$ into $L(Y,P^\pi)$ by showing that it preserves edges.
     
     Let $Q \subseteq X^d$ be a $d$-dimensional $k$-template that is an edge of $L(X^d,P)$. Let $g : P \into Q$ be a bijection demonstrating that $Q$ is a homomorphic image of $P$; that is, if $x,y \in P$ and $x_i = y_i$, then $g(x)_i = g(y)_i$. Define $h : P^\pi \into Y$ so that $h(y) = f_\pi gf_\pi^{-1}(y)$. Then, $h[P^\pi] = f_\pi[Q]$. It is  easily seen that $h$ demonstrates that $f_\pi[Q]$ is a homomorphic image of $P$. Thus, $Q$ is an edge of $L(Y,P^\pi)$. \qed
     
     \bigskip

    {\sc Corollary 4.3}: {\em Let $H$ be an algebraic $k$-hypergraph and $P$ be a $d$-dimensional 
    $k$-template. The following are equivalent$:$
    
    $(1)$ For every $M < \omega$, $L(M^d,P)$ is embeddable into $H$.
    
    $(2)$ $L(\RR^d,P)$ is embeddable into $H$.
    
    $(3)$ $L(\RR^d,P)$ is immersible into $H$.}

    \bigskip
    
    {\it Proof}.  We will prove 
    $(1) \Longrightarrow (2) \Longrightarrow(3) \Longrightarrow (1)$.

    \smallskip
    
    $(1) \Longrightarrow (2)$: 
    Assume $(1)$ is true. 
    By Theorem~3.1, there are $\pi : d \into m$ and  a semialgebraic immersion  of 
    $L(\RR^m,P^\pi)$ into $H$. Since $L(\RR^m,P^\pi)$ is semialgebracally immersible into $H$, then by Lemma~2.3, it is also embeddable into $H$. But Lemma~4.2 implies that 
    $L(\RR^d,P)$ is embeddable into $L(\RR^m,P^\pi)$, so $L(\RR^d,P)$ is embeddable into $H$. 
    Thus, $(2)$ holds. 
    
    \smallskip
    
    $(2) \Longrightarrow (3)$: This is trivial since any embedding of $L(\RR^d,P)$  into $H$ is also an immersion.
    
    \smallskip
    
    $(3) \Longrightarrow (1)$: Assume $(3)$ and let $f$ be an immersion of $L(\RR^d,P)$  into $H$.
    Then $f$ is one-to-one in each coordinate. Consider $M < \omega$. Lemma~4.1 implies that there 
    are $Y_0,Y_1, \ldots, Y_{d-1} \in [\RR]^M$ such that $f$ is one-to-one on $Y = Y_0 \times Y_1 \times \cdots \times Y_{d-1}$. Then, $f \harp Y$ is an embedding of $L(Y,P)$ into $H$.
    Since $L(M^d,P) \cong L(Y,P)$, then $L(M^d,P)$ also  is embeddable into $H$
    \qed
    
        \bigskip

 {\sc Definition 4.4}: If $H$ is any $k$-hypergraph, then we define $\dep(H)$, the  {\bf depth} of $H$,  to be the least $\delta < \omega$ 
 for which there is a $d$-dimensional $k$-template $P$ such that $e(P) = \delta +1$ and 
 $H$ embeds $L(M^d,P)$  for all $M < \omega$. If there is no such $P$, then  $\dep(H)= {\infty}$. 
 
 \bigskip
 
 {\sc Remark:} If $P$ is a $d$-dimensional $k$-template, then, according to Lemma~2.2,  there is an $e(P)$-dimensional $k$-template $Q$ such that for every set $X$, $L(X^d,P)$ embeds 
 $L(X^{e(P)},Q)$. Therefore, we get the following alternative definition of depth: {\em If $H$ is a $k$-hypergraph and $\dep(H) = \delta < \omega$,  then $\delta$ is  the least  
 for which there is a $(\delta+1)$-dimensional $k$-template $P$ such that  $H$ embeds $L(M^{\delta+1},P)$  for all $M < \omega$.}
 
 \bigskip
 
 In the case of algebraic hypergraphs, depth can be defined by substituting the single hypergraph $L(\RR^d,P)$ for all the finite $L(M^d,P)$. (See Corollary~4.3.)
 
 \bigskip
 
 {\sc Corollary 4.5}: {\em Suppose that $H$ be an algebraic $k$-hypergraph and $\delta = \dep(H)$. If $\delta < \omega$, then  
 $\delta$ is the least for which there is a $d$-dimensional $k$-template $P$ such that $e(P) = \delta +1$ and $H$ embeds 
 $L(\RR^d,P)$. If $\delta = \infty$, then there  is no $P$ 
 such that  $H$ embeds $L(\RR^d,P)$.}
 \bigskip
 
 {\it Proof}. Suppose $\delta < \omega$. First, suppose that $P$ is a $d$-dimensional $k$-template such that $H$ embeds $L(\RR^d,P)$. Then $L(M^d,P)$ is embeddable into $H$ for every $M < \omega$. Thus, $e(P) \geq \delta +1$. Next, let $P$ be a $d$-dimensional $k$-template 
 such that $L(M^d,P)$ is embeddable into $H$ for all $M < \omega$. By Theorem~3.1, there is a surjection $\pi : d \into m$ such that  $L(\RR^m,P^\pi)$ is semialgebraically immersible into $H$.
 By Corollary~4.3,  $L(M^m,P^\pi)$ is embeddable into $H$ for every $M < \omega$. Then, $e(P^\pi) \leq e(P)$. On the other hand, by the minimality of $\delta$, $e(P^\pi) \geq e(P)$. Thus $e(P^\pi) = \delta +1$ and $L(\RR^m,P^\pi)$ is embeddable into $H$. 
 
 \smallskip
 
 Suppose $\delta = \infty$. If  $P$ were such that $L(\RR^d,P)$ is embeddable into $H$, then, 
 for all $M < \omega$,  
 $L(M^d,P)$  would be embeddable into $H$, and then $\delta < \omega$. \qed
 
 \bigskip

 The notion of depth was suggested by the following characterization of infinite chromatic numbers of 
 algebraic hypergraphs. 
 
\bigskip
 
 {\sc Corollary 4.6}: {\em Suppose that $H$ is an algebraic $k$-hypergraph, $\delta = \dep(H)$  and 
 $\kappa$  is an infinite cardinal. 
Then, $\chi(H) \leq \kappa$ iff $\kappa^{+\delta} \geq 2^{\aleph_0}$.} 

\bigskip

{\it Proof}. This is an immediate consequence of Corollary~4.3 and Lemmas~ 2.1 and~2.2. \qed

\bigskip

We next extend the definition of depth to polynomials: 
 if $R$ is a real closed  field and $p(x_0,x_1, \ldots, x_{k-1})$ is a $(k,n)$-ary polynomial over $R$, then  the {\bf depth} of $p(x_0,x_1, \ldots,x_{k-1})$ is the depth of its zero $k$-hypergraph in $R$. 
 
 Let ${\mathcal Q}$ be the set of polynomials over $\QQ$ that are $(k,n)$-ary for some $k,n < \omega$. Clearly, ${\mathcal Q}$ is computable. 
 
 \bigskip
 
 {\sc Corollary 4.7}: {\em The depth function restricted to ${\mathcal Q}$ is computable.} 
 
 \bigskip
 
 {\it Proof}. This proof implicitly makes use of the decidability of ${\mathsf {RCF}}$. 
 
 Suppose we are given a $(k,n)$-ary polynomial $p(x_0,x_1, \ldots, x_{k-1})$ and its 
 zero $k$-hypergraph $H$. Of course, the depth of $p(x_0,x_1, \ldots, x_{k-1})$is $\dep(H)$. By the Remark following Definition~4.4, we need only be concerned with the finitely many $P$ that are $d$-dimensional $k$-templates with $d < k$. For each such $P$,  do an effective  search to get either (1) an $M < \omega$ such that $H$ does not embed $L(M^d,P)$ or else (2) a surjection 
 $\pi : d \into m$ and 
 an ${\mathcal L}_{\sf OF}$-formula that defines an immersion of $L(\RR^m,P^\pi)$ into $H$. 
 By Theorem~3.1, this search will terminate. If, for every such $P$, the search terminates with (1),  then $\dep(H) = \infty$. Otherwise, let $d$ be the least for which there is some 
 $d$-dimensional $k$-template $P$ for which the search terminates in (2). Then, $\dep(H)= d-1$. 
 \qed
  
 \bigskip
 
  The next corollary becomes the conjecture from~\cite{avoid}  in the special case that $\kappa = \aleph_0$,

\bigskip
 
 {\sc Corollary 4.8}: {\em For each infinite cardinal $\kappa$, the set of $\kappa$-avoidable polynomials in ${\mathcal Q}$ is computable.}

    \bigskip
    
    {\it Proof}. Fix infinite $\kappa$ and let $m< \omega$ be the least such that $\kappa^{+m} \geq 2^{\aleph_0}$ if possible. Otherwise, let $m = \infty$. 
    
    For a given polynomial in ${\mathcal Q}$, let $\delta$ be its depth. Then, by Corollary~4.6,  the polynomial is $\kappa$-avoidable iff $\delta \geq m$. \qed

   \bigskip
   
   The instance of Theorem~3.1 in which $d=1$ is an interesting special case. 
 If $A$ is a set and $k < \omega$, then $K^{(k)}_A$ is the complete $k$-hypergraph on $A$; 
 that is, $K^{(k)}_A = (A,[A]^k)$. 
 In particular, $K^{(2)}_A = K_A$ is the complete graph on $A$.

 \bigskip
 
 {\sc Corollary~4.9}: {\em If  $H = (\RR^n,E)$ is an algebraic $k$-hypergraph, then 
 the following are equivalent$:$
 
 $(1)$ for each $m < \omega$, $K^{(k)}_m$ is embeddable into $H;$
 
 $(2)$ there is an infinite algebraic $A \subseteq \RR^n$ such that $K^{(k)}_A$ is a subhypergraph of $H.$}  \qed
 
 \bigskip

 {\sc Remark}: Being infinite and algebraic, the set $A$ in $(2)$ has cardinality $2^{\aleph_0}$. 
 
 \bigskip
 
 {\it Proof}. Obviously, $(2) \Longrightarrow (1)$.
 
 For the converse, suppose that $(1)$ holds. Let $p(x_0,x_1, \ldots,x_{k-1})$ be 
 a symmetric, reflexive $(k,n)$-ary polynomial whose zero hypergraph is $H$. 
 Let $P$ be a $1$-dimensional $k$-template. 
 (Note that, up to isomorphism, there is only one.) For every set $X$, $L(X,P) \cong K^{(k)}_X$. 
 Thus, Theorem~3.1 implies that  
 $L(\RR,P)$ is semialgebraically immersible into $H$. Therefore,  let $f : \RR \into  \RR^n$ be a Nash immersion of $L(\RR,P)$ into $H$, which necessarily is an embedding. 
 Let $B$ be the range of $f$. Then, whenever $b_0,b_1, \ldots, b_{k-1} \in B$, then $p(b_0,b_1, \ldots, b_{k-1}) = 0$. By the Hilbert Basis Theorem, let $A \subseteq \RR^n$ be a minimal algebraic subset of $\RR^n$  such that $A \cap B$ is infinite.  
 We will show that whenever $a_0,a_1, \ldots, a_{k-1} \in A$, then $p(a_0,a_1, \ldots,a_{k-1}) = 0$.
 
 The proof of the claim is by induction. For $c \leq k$, let ${\sf S}(c)$ be the statement:
 
 \begin{quote}
 
 Whenever  $a_0,a_1, \ldots, a_{k-1} \in A$ and 
\mbox{$|\{i < k : a_i \not\in B\}| \leq c$}, then $p(a_0,a_1, \ldots,a_{k-1}) = 0$.

\end{quote}
Our goal is to prove ${\sf S}(k)$. Clearly, ${\sf S}(0)$ is true.

Suppose that $c < k$ and that ${\sf S}(c)$ is true. Suppose that $a_0,a_1, \ldots, $ $a_{k-1} \in A$ and that 
$a_0,a_1, \ldots, a_c \not\in B$. Consider the set $X = \{x \in A : p(x,a_1,a_2, \ldots,a_{k-1}) = 0\}$. 
By the inductive hypothesis, $X$ is infinite, so $A \subseteq X$. Therefore, $a_0 \in X$ and 
$ p(a_0,a_1,a_2, \ldots,a_{k-1}) = 0\}$. Thus, ${\sf S}(c+1)$. \qed

\bigskip

{\sc Remark}: If an algebraic $k$-hypergraph $H$ satisfies $(1)$ and/or $(2)$ of Corollary~4.9, then $\chi(H) = 2^{\aleph_0}$. The converse does not hold if $2^{\aleph_0}$ is a limit cardinal, for then 
$\chi(H) = 2^{\aleph_0}$ iff $\chi(H)$ is uncountable.  However, if $2^{\aleph_0}$ is a successor cardinal and $\chi(H) = 2^{\aleph_0}$, then  $(1)$ and/or $(2)$. 

\bigskip

In order to generalize Corollary~4.7 from ${\mathcal Q}$ to all polynomials over $\RR$, we make the following definition.  \bigskip
 
 {\sc Definition 4.10}:
 A collection ${\mathcal H}$ of algebraic hypergraphs is {\bf decidable} if there is a computable function 
that maps each $((k,n)+\ell)$-ary polynomial $p(x_0,x_1, \ldots,x_{k-1},y)$ over $\QQ$ to an $\ell$-ary formula $\theta(y)$ such that 
if $c \in \RR^\ell$ and $H_c$ is the zero $k$-hypergraph of $p(x_0,x_1, \ldots,x_{k-1},c)$, 
the $H_c \in {\mathcal H}$ iff  ${ \RR} \models \theta(c)$.

 \bigskip

 {\sc Corollary 4.11}: {\em For each $\delta \in \omega \cup \{\infty\}$, the set of 
 algebraic hypergraphs having depth $\delta$ is decidable uniformly in $\delta$.}
 
 \bigskip
 
 {\it Proof}. Suppose that we are given a $((k,n)+\ell)$-ary polynomial in ${\mathcal Q}$. 
 Let $H_y$ be its zero $k$-hypergraph. 
 For each $d$-dimensional $k$-template $P$ and surjection $\pi : d \into m$, where $m \leq d < k$, effectively get $M$  as in Theorem~3.1. We can get the same $M$ for every $P$.  If $\delta = \infty$, let $\theta(y)$ be the formula asserting: for no $P$ does $H_y$ embed $L(M^d,P)$. If $\delta < \omega$, then let $\theta(y)$ assert: $\delta$ is the least for which there is   $P$ $e(P) = \delta+1$ and 
 $H_y$ embeds $L(M^d,P)$. \qed 
  
 \bigskip
 
 {\sc Corollary 4.12}: {\em If $\kappa$ is an infinite cardinal, then the set of algebraic $\kappa$-colorable hypergraphs is decidable.}
 
 \bigskip
 
 {\it Proof.} This is immediate from Corollaries~4.6 and 4.11. \qed
 
 \bigskip

{\sc Corollary 4.13}: {\em Let $\kappa$ be an uncountable cardinal. The set of algebraic hypergraphs $H$ such that $\chi(H) = \kappa$ is decidable.} \qed

 \bigskip
 
 We end  with the  following open question suggested by the previous corollary. 
 
 \bigskip
 
 {\sc Question 4.14}: Is the collection of algebraic hypergraphs $H$ such that $\chi(H) = \aleph_0$ decidable?

  \bigskip
 
 We know that the set of those polynomials over $\QQ$ having a zero hypergraphs with chromatic number $\aleph_0$ is ${\Pi_2^0}$. Can that be improved?
 
  \bigskip

\bibliographystyle{plain}

\end{document}